\documentclass[11pt,amsfonts]{article}
\usepackage{graphicx}
\usepackage{latexsym}
\usepackage{amssymb}
\usepackage{amsmath}
\usepackage{layout}
\newtheorem{prop}{Proposition}
\newtheorem{lemma}{Lemma}
\newtheorem{definition}{Definition}

\newtheorem{theorem}{Theorem}
\newtheorem{remark}{Remark}

\def\real{{\mathord{{\rm I\kern-2.8pt R}}}}        
\def\inte{{\mathord{{\rm I\kern-2.8pt N}}}}

\def\sZZ{{\rm Z\kern-2.8ptem{}Z}}

\def\z{{\mathchoice
  {\sZZ}
  {\sZZ}
  {\rm Z\kern-0.30em{}Z}
  {\rm Z\kern-0.25em{}Z} }}
\def\sQQ{{\kern 0.27em \vrule height1.45ex width0.03em depth0em
          \kern-0.30em \rm Q}}
\def\qu{{\mathchoice
    {\sQQ}
    {\sQQ}
  {\kern 0.225em \vrule height1.05ex width0.025em depth0em \kern-0.25em \rm Q}
  {\kern 0.180em \vrule height0.78ex width0.020em depth0em \kern-0.20em \rm Q}
        }}
\def\sCC{{\kern 0.27em \vrule height1.45ex width0.03em depth0em
          \kern-0.30em \rm C}}
\def\complex{{\mathchoice
    {\sCC}
    {\sCC}
  {\kern 0.225em \vrule height1.05ex width0.025em depth0em \kern-0.25em \rm C}
  {\kern 0.180em \vrule height0.78ex width0.020em depth0em \kern-0.20em \rm C}
        }}

\newcommand{\ba}{\begin{array}}
\newcommand{\ea}{\end{array}}
\newcommand{\be}{\begin{equation}}
\newcommand{\ee}{\end{equation}}
\newcommand{\bea}{\begin{eqnarray}}
\newcommand{\eea}{\end{eqnarray}}
\newcommand{\beaa}{\begin{eqnarray*}}
\newcommand{\eeaa}{\end{eqnarray*}}

\def\b{\beta}

\def\z{\zeta}

\font\tenmath=msbm10 \font\sevenmath=msbm7 \font\fivemath=msbm5
\newfam\mathfam \textfont\mathfam=\tenmath
\scriptfont\mathfam=\sevenmath \scriptscriptfont\mathfam=\fivemath

\def \b{\noindent}

\def \={{\buildrel {\rm (law)} \over =}}

\def\qed{ \hfill \vrule width.25cm height.25cm depth0cm\smallskip}

\newcommand{\basa}{\begin{assumption}}
\newcommand{\easa}{\end{assumption}}

\newcommand{\bas}{\begin{assum}}
\newcommand{\eas}{\end{assum}}

\newcommand{\ignore}[1]{}
\begin{document}

\renewcommand{\thefootnote}{\fnsymbol{footnote}}

\date{ }

\title{2D- stochastic currents over the Wiener  sheet \footnote{Dedicated to the memory of Constantin Tudor. }}
\author{Franco Flandoli$^{1}, $\qquad  Peter Imkeller  $^{2}\qquad $%
Ciprian A. Tudor $^{3, 4}$ \footnote{ Supported by the CNCS grant PN-II-ID-PCCE-2011-2-0015. Associate member of the team Samos, Universit\'e de Panth\'eon-Sorbonne Paris 1.  The author would like to acknowledge generous support from the Alexander von Humboldt Foundation which made possible several research visits at the Humboldt Universit\"at zu Berlin. The support by  the ANR grant "Masterie" BLAN 012103 is also acknowledged.  }
\vspace*{0.1in} \\
$^{1}$Dipartimento di Matematica Applicata, Universita di Pisa,\\ Via Bonnano 25B, I-56126, Pisa, Italy
\\flandoli@dma.unipi.it\vspace*{0.1in}\\
$^{2}$Institut f\"ur Mathematik, Humboldt-Universit\"at zu Berlin\\
, Unter den Linden 6, 10099 Berlin, Germany,\\
 imkeller@mathematik.hu-berlin.de \vspace*{0.1in} \\
 $^{3}$ Laboratoire Paul Painlev\'e, Universit\'e de Lille 1\\
 F-59655 Villeneuve d'Ascq, France,\vspace*{0.1in} \\
  $^{4}$ Department of Mathematics, \\
Academy of Economical Studies, Bucharest, Romania \vspace*{0.1in}\\
 \quad tudor@math.univ-lille1.fr\vspace*{0.1in}}
\maketitle

\begin{abstract}
By using stochastic calculus for two-parameter processes and chaos expansion into multiple Wiener-It\^o integrals, we define a 2D-stochastic current over the Brownian sheet. This concept comes from geometric measure theory. We also study the regularity of the stochastic current with respect to the randomness variable in the Watanabe spaces and with respect to the spatial variable in the deterministic Sobolev spaces.
\end{abstract}

\vskip0.5cm

{\bf  2000 AMS Classification Numbers:} 60G15, 60G18,  60H05, 76M35, 60H05

 \vskip0.3cm


\section{Introduction}

The purpose of this paper is to study two-dimensional stochastic currents. The concept of current comes from geometric measure theory.
Let us  briefly recall a few definitions on deterministic currents used in
the sequel. More informations can be found in \cite{Fed}, \cite{Morgan}, \cite{Simon},
\cite{GMS}.

We shall denote the Euclidean norm and scalar product in $\mathbb{R}^{d}$ ($d$
is fixed throughout the paper) by $\left|  \cdot\right|  $ and $\left\langle
\cdot,\cdot\right\rangle $ respectively. We shall also need the space
$\Lambda_{2}\mathbb{R}^{d}$ of $2$-vectors, its dual space $\Lambda
^{2}\mathbb{R}^{d}$ of $2$-covectors, and the duality between them, still
denoted by $\left\langle \cdot,\cdot\right\rangle $. If we represent vectors
and covectors in the standard basis (the summations are extented from $1$ to
$d$)
\[
v=\sum_{i<j}v^{ij}\,e_{i}\wedge e_{j},\qquad w=\sum_{i<j}w_{ij}\,e^{i}\wedge
e^{j}%
\]
then
\[
\left\langle w,v\right\rangle =\sum_{i<j}w_{ij}v^{ij}.
\]
Notationally, we prefer to write the covectors in the first argument of
$\left\langle .,.\right\rangle $ since later on the notation for stochastic
integrals is more natural. The norm of $2$-vectors and $2$-covectors is also
defined as
\[
\left|  v\right|  ^{2}=\sum_{i<j}\left(  v^{ij}\right)  ^{2},\qquad\left|
w\right|  ^{2}=\sum_{i<j}\left(  w_{ij}\right)  ^{2}.
\]
Notice that for all $v_{1},v_{2}\in\mathbb{R}^{d}$ and $w^{1},w^{2}%
\in\mathbb{R}^{d}$, the duality between the $2$-vector
\[
v_{1}\wedge v_{2}=\sum_{i<j}\left[  \left(  v_{1}\right)  _{i}\left(
v_{2}\right)  _{j}-\left(  v_{1}\right)  _{j}\left(  v_{2}\right)
_{i}\right]  \,e_{i}\wedge e_{j}%
\]
and the $2$-covector $w^{1}\wedge w^{2}$ similarly defined, is given by
\[
\left\langle w^{1}\wedge w^{2},v_{1}\wedge v_{2}\right\rangle =\det\left(
\left\langle v_{i},w^{j}\right\rangle \right)  .
\]
Similarly we have
\[
\left|  v_{1}\wedge v_{2}\right|  ^{2}=\det\left(  \left\langle v_{i}%
,v_{j}\right\rangle \right)  =\left|  v_{1}\right|  ^{2}\left|  v_{2}\right|
^{2}-\left\langle v_{1},v_{2}\right\rangle ^{2}.
\]
Let $\mathcal{D}^{k}$ be the space of all infinitely differentiable and
compactly supported $k$-forms on $\mathbb{R}^{d}$. A $k$\textit{-dimensional
current} is a linear continuous functional on $\mathcal{D}^{k}$. We denote by
$\mathcal{D}_{k}$ the space of $k$-currents. In this paper we are only
interested in the cases $k=1,2$, so we recall a few corresponding notations.
Let $\varphi$ (resp. $\psi$) be an element of $\mathcal{D}^{1}$ (resp.
$\mathcal{D}^{2}$). We shall write them as
\[
\varphi=\sum_{i=1}^{d}\varphi_{i}\,dx^{i},\qquad\psi=\sum_{i<j}\psi
_{ij}\,dx^{i}\wedge dx^{j}%
\]
where $\varphi_{i}$ and $\psi_{ij}$ are infinitely differentiable and
compactly supported functions on $\mathbb{R}^{d}$. Typical examples of
$1$-currents are those induced by regular curves $\left(  X_{t}\right)
_{t\in\left[  0,T\right]  }$ in $\mathbb{R}^{d}$:
\begin{equation}\label{map1}
\varphi\mapsto S(\varphi):=\int_{0}^{T}\left\langle \varphi\left(
X_{t}\right)  ,\dot{X}_{t}\right\rangle \,dt
\end{equation}
while typical examples of $2$-currents are given by regular surfaces $\left(
f\left(  t,s\right)  \right)  _{(t,s)\in A}$ in $\mathbb{R}^{d}$, where $A$ is
a Borel set in $\mathbb{R}^{2}$:
\begin{eqnarray}\label{map2}
\psi &&  \mapsto T(\psi):=\int_{A}\left\langle \psi\circ f,\frac{\partial
f}{\partial t}\wedge\frac{\partial f}{\partial s}\right\rangle \,dtds\nonumber \\
&&  =\int_{A}\sum_{i<j}\left(  \psi_{ij}\circ f\right)  \,\left\langle
dx^{i}\wedge dx^{j},\frac{\partial f}{\partial t}\wedge\frac{\partial
f}{\partial s}\right\rangle \,dtds \nonumber \\
&&  =\int_{A}\sum_{i<j}\left(  \psi_{ij}\circ f\right)  \,\left(
\frac{\partial f_{i}}{\partial t}\frac{\partial f_{j}}{\partial s}%
-\frac{\partial f_{j}}{\partial t}\frac{\partial f_{i}}{\partial s}\right)
\,dtds.
\end{eqnarray}

The generalization of the previous example of $1$-current to the stochastic
case has been the object of previous research in the recent  literature. We refer to \cite{Gia}, \cite{FG}, \cite{FGR}, \cite{FT} for the study of 1D-stochastic currents over the Brownian motion or over the fractional Brownian motion using various interpretation (It\^o, Skorohod, Stratonovich etc) for the stochastic integral appearing in the definition of the stochastic current.
  We will investigate    in our paper the 2D-current   in the stochastic context using the theory of stochastic integration  for two -parameter stochastic processes which has been developed in the eighties (see e.g. \cite{Dozzi}, \cite{Imk}, \cite{N81}). Somehow, our analysis related geometric aspects of stochastic processes with Malliavin calculus in the spirit of \cite{Mall}.

Let $(W_{s,t} )_{\left(  t,s\right)
\in\left[  0,1\right]  ^{2}}$ be a $\left(  2,d\right)  $-Brownian sheet, with
$d>2$. By this we mean a centered Gaussian $d$-dimensional process with
$2$-dimensional time parameter $\left(  t,s\right)  \in\left[  0,1\right]  ^{2}$
such that
\[
Cov\left(  W^{i}_{s,t}, W^{j}_{s^{\prime}, t^{\prime}} \right)  =\delta_{ij}\left(  s\wedge t\right)  \left(  s^{\prime
}\wedge t^{\prime}\right)
\]
where $\delta_{ij}$ denotes the Kronecker symbol.
Note that the functional (\ref{map1}) defines a vector valued distribution
$$\int_{0} ^{1} \Delta (x -X_{t})dX_{t}$$
with $\Delta$ the Dirac distribution and the meaning of $ \Delta (x -X_{t})$ will be clarified  later.
Given (\ref{map2}) the correct multiparameter stochastic integrals to be studied for the theory
of currents are formally written as
\[
T_{ij}\left(  x\right)  :=\int_{\left[  0,1\right]  ^{2}}\Delta\left(
x-W_{s,t}  \right)  \partial_{t}W^{i}_{s,t}
\partial_{s}W^{j}_{s,t}
\]
with $1\leq i\neq j\leq d$. An important question is how to define the above stochastic integral.

At this level, there is no restriction to take $i=1$, $j=2$. Thus we need to define and to study
study the mapping
\begin{equation}\label{int}
T\left(  x\right)  :=\int_{\left[  0,1\right]  ^{2}}\Delta\left(  x-W_{s,t}  \right)  \partial_{t}W^{1}_{s,t} \partial_{s}%
W^{2}_{s,t}.
\end{equation}
This is actually the aim of this work: to give a meaning to the integral (\ref{int}) and to study its regularity  in the  deterministic Sobolev spaces $H^{-r} (\mathbb{R} ^{d}; \mathbb{R} ^{d})$ as a function of $x$ and in the Sobolev-Watanabe spaces $\mathbb{D}^{\alpha , 2}$ as a function of $\omega$. As we mentioned above, the idea to define the stochastic integral (\ref{int}) comes from the stochastic integration theory for two-parameter processes. We will explain it in details in Section \ref{defint}. This integral also involves the divergence (Skorohod) integral with respect to the Brownian sheet and multiple Wiener-It\^o integrals. Therefore, before defining (\ref{int}) we will recall in the next section some elements of the Malliavin calculus for Gaussian processes. We will find the following answer: although the definition of the 2D -stochastic current is substantially different from the 1D-case, the 2D-current inherits the regularity of the one-dimensional case, with respect to $\omega$ as well as  with respect to
$x$. That means the following: the  integral (\ref{int}) belongs to the Sobolev -Watanabe space $\mathbb{D} ^{-\alpha , 2} $ for every $\alpha >\frac{1}{2}$ and to the Sobolev space $H^{-r}(\mathbb{R}; \mathbb{R})$ for $r> \frac{1}{2}$. We notice the same interesting phenomenon already remark in the one-dimensional case (see \cite{FT}): the order of regularity of (\ref{int}) is the same with respect to $\omega$ and with respect to $x$.

 Our paper is structured as follows: Section 2 contains some preliminaries on multiple Wiener-It\^o integrals and the divergence integral.  In Section 3 we explain how we arrive to the our definition of the 2D -stochastic currents, making links with the stochastic calculus for two-parameter processes. In Sections 4 and 5 we analyze the 2D-stochastic current over the Brownian sheet in Watanabe spaces and Sobolev spaces respectively. Finally in Section 6 we treat the 2D -stochastic current over the $d$ -dimensional Brownian sheet.

\section{Malliavin calculus, chaos expansion and Sobolev-Watanabe spaces}\label{prel}

Here we describe the elements from stochastic analysis that we will
need in the paper. Consider ${\mathcal{H}}$ a real separable Hilbert
space and $(B (\varphi), \varphi\in{\mathcal{H}})$ an isonormal
Gaussian process on a probability space $(\Omega, {\cal{A}}, P)$, that is a centered Gaussian family of random
variables such that $\mathbf{E}\left( B(\varphi) B(\psi) \right)  =
\langle\varphi, \psi\rangle_{{\mathcal{H}}}$.
Denote by $I_{n}$ the multiple stochastic integral with respect to
$B$ (see \cite{N}). This $I_{n}$ is actually an isometry between the
Hilbert space ${\mathcal{H}}^{\odot n}$(symmetric tensor product)
equipped with the scaled norm
$\frac{1}{\sqrt{n!}}\Vert\cdot\Vert_{{\mathcal{H}}^{\otimes n}}$ and
the Wiener chaos of order $n$ which is defined as the closed linear
span of the random variables $H_{n}(B(\varphi))$ where
$\varphi\in{\mathcal{H}}, \Vert\varphi\Vert_{{\mathcal{H}}}=1$ and
$H_{n}$ is the Hermite polynomial of degree $n\geq 1$
\begin{equation*}
H_{n}(x)=\frac{(-1)^{n}}{n!} \exp \left( \frac{x^{2}}{2} \right)
\frac{d^{n}}{dx^{n}}\left( \exp \left( -\frac{x^{2}}{2}\right)
\right), \hskip0.5cm x\in \mathbb{R}.
\end{equation*}
The isometry of multiple integrals can be written as: for $m,n$ positive integers,
\begin{eqnarray}
\mathbf{E}\left(I_{n}(f) I_{m}(g) \right) &=& n! \langle f,g\rangle _{{\mathcal{H}}^{\otimes n}}\quad \mbox{if } m=n,\nonumber \\
\mathbf{E}\left(I_{n}(f) I_{m}(g) \right) &= & 0\quad \mbox{if } m\not=n.\label{iso}
\end{eqnarray}
It also holds that
\begin{equation*}
I_{n}(f) = I_{n}\big( \tilde{f}\big)
\end{equation*}
where $\tilde{f} $ denotes the symmetrization of $f$ defined by $\tilde{f}%
(x_{1}, \ldots , x_{x}) =\frac{1}{n!} \sum_{\sigma \in {\cal S}_{n}}
f(x_{\sigma (1) }, \ldots , x_{\sigma (n) } ) $.

We recall that any square integrable random variable which is
measurable with respect to the $\sigma$-algebra generated by $B$ can
be expanded into an orthogonal sum of multiple stochastic integrals
\begin{equation}
\label{sum1} F=\sum_{n\geq0}I_{n}(f_{n})
\end{equation}
where $f_{n}\in{\mathcal{H}}^{\odot n}$ are (uniquely determined)
symmetric functions and $I_{0}(f_{0})=\mathbf{E}\left[  F\right]  $.

Let $L$ be the Ornstein-Uhlenbeck operator
\begin{equation*}
LF=-\sum_{n\geq 0} nI_{n}(f_{n})
\end{equation*}
if $F$ is given by (\ref{sum1}).

For $p>1$ and $\alpha \in \mathbb{R}$ we introduce the
Sobolev-Watanabe space $\mathbb{D}^{\alpha ,p }$  as the closure of
the set of polynomial random variables with respect to the norm
\begin{equation*}
\Vert F\Vert _{\alpha , p} =\Vert (I -L) ^{\frac{\alpha }{2}} \Vert
_{L^{p} (\Omega )}
\end{equation*}
where $I$ represents the identity. In this way, a random variable
$F$  as in (\ref{sum1}) belongs $\mathbb{D}^{\alpha , 2}$ if and
only if
\begin{equation}\label{d12}
\sum_{n\geq 0} (1+n) ^{\alpha } \Vert I_{n}(f_{n}) \Vert
_{L^{2}(\Omega)} ^{2} =\sum_{n\geq 0} (1+n) ^{\alpha }n! \Vert f_{n}
\Vert ^{2} _{{\cal{H}}^{\otimes n}}<\infty .
\end{equation}
We denote by $D$  the Malliavin  derivative operator that acts on smooth functions of the form $F=g(B(\varphi _{1}), \ldots , B(\varphi_{n}))$ ($g$ is a smooth function with compact support and $\varphi_{i} \in {{\cal{H}}}$)
\begin{equation*}
DF=\sum_{i=1}^{n}\frac{\partial g}{\partial x_{i}}(B(\varphi _{1}), \ldots , B(\varphi_{n}))\varphi_{i}.
\end{equation*}
The operator $D$ is continuous from $\mathbb{D} ^{\alpha , p} $ into
$\mathbb{D} ^{\alpha , p} \left( {\cal{H}}\right).$ The adjoint of
$D$ is denoted by $\delta $ and is called the divergence (or
Skorohod) integral. It is a continuous operator from $\mathbb{D}
^{\alpha, p } \left( {\cal{H}}\right)$ into $\mathbb{D} ^{\alpha ,
p}$. For adapted integrands, the divergence integral coincides to
the classical It\^o integral. We will use the notation
\begin{equation*}
\delta (u) =\int_{0}^{T} u_{s} dB_{s}.
\end{equation*}
Let  $u$ be  a stochastic process having the chaotic decomposition $u_{s}=\sum _{n\geq 0} I_{n}(f_{n}(\cdot ,s))$
where $f_{n}(\cdot, s)\in {\cal{H}}^{\otimes n}$ for every $s$. One
can prove that $u \in {\rm Dom} \ \delta$ if and only if $\tilde f_n
\in {\cal{H}}^{\otimes (n+1)}$ for every $n \geq 0$, and
$\sum_{n=0}^{\infty}I_{n+1}(\tilde f_n)$ converges in $L^2(\Omega)$.
In this case,
\begin{equation}
\label{sko}\delta(u)=\sum_{n=0}^{\infty}I_{n+1}(\tilde f_n) \quad \mbox{and}
\quad \mathbf{E}|\delta(u)|^2=\sum_{n=0}^{\infty}(n+1)! \ \|\tilde
f_n\|_{{\cal{H}}^{\otimes (n+1)}}^{2}.
\end{equation}
We will also recall the integration by parts formula
\begin{equation}
\label{ip}
F\delta (u)= \delta (Fu)+ \langle DF, u\rangle _{{\cal{H}}}
\end{equation}
which holds for $F\in \mathbb{D} ^{\alpha, 2}, u\in \mathbb{D} ^{\alpha , 2} ({\cal{H}}) $ and $Fu \in \mathbb{D}^{\alpha ,2} ({\cal{H}})$.
In the present work we will consider divergence integral with respect to a Brownian motion in $\mathbb{R}^{d}$.
Throughout this paper we will denote  by $p_{s}(x) $ the Gaussian
kernel of variance $s>0$ given by $p_{s}(x)= \frac{1}{\sqrt{2\pi s}} e^{-\frac{x^{2}}{2s}},
x\in \mathbb{R}$ and for $x=(x_{1}, \ldots , x_{d})\in \mathbb{R}^{d}$ by $p_{s} ^{d}(x)= \prod_{i=1}^{d}p_{s}(x_{i}).$

\section{Definition of the
 2D- stochastic current}\label{defint}

The idea to define a stochastic integral  (\ref{int}) comes from the stochastic calculus for two-parameter martingales developed in the eighties. Let us briefly recall the context. Consider a two- parameter Wiener process $(W_{s,t})_{s,t\in [0,1]}$ and  let  us express the random variable $W_{s,t}^{2}$ (here $W_{s,t}$ is one dimensional random variable and $W_{s,t}^{2}$ is its square; we mention this in order to avoid the confusion with the components of the two-dimensional Brownian sheet that appears in Section \ref{prel})  as a stochastic integral plus a deterministic integral in the spirit of the classical It\^o formula.  Let $t,s \in [0,1] $ be fixed  and let
\begin{equation*}
0=t_{0}<t_{1}<\ldots t_{n}=t
\mbox{ and }
0=s_{0}<s_{1}<\ldots s_{m} = s
\end{equation*}
be two partitions of the intervals $[0,t]$ and $[0,s]$ respectively. For fixed $t$, we can write
\begin{eqnarray*}
W_{s,t}^{2}= 2\sum_{i=0}^{m-1} W_{s_{i}, t }\left(W_{s_{i+1}, t }-W_{s_{i}, t }\right)+ \sum_{i=0}^{m-1} \left( W_{s_{i+1}, t }-W_{s_{i}, t }\right) ^{2}  .
\end{eqnarray*}
It is standard to prove that the last summand $\sum_{i=0}^{m-1} \left( W_{s_{i+1}, t }-W_{s_{i}, t }\right) ^{2}  $ converges in $L^{2}(\Omega)$ as $m\to \infty$ to
$st$ which, in some sense, constitutes the quadratic variation of the two-parameter martingale $W_{s,t}$.  Let us analyze the first summand. Writing $g(t)= g(0)+ \sum_{j=0}^{n-1} (g(t_{j+1}) -g(t_{j}))$ we get
\begin{eqnarray*}
\sum_{i=0}^{m-1} W_{s_{i}, t }\left(W_{s_{i+1}, t }-W_{s_{i}, t }\right)
&=& \sum_{i=0}^{m-1} \sum_{j=0} ^{m-1} W_{s_{i}, t_{j+1} } \left( W_{s_{i+1}, t_{j+1} }-W_{s_{i+1}, t_{j}}-W_{s_{i}, t_{j+1} } +W_{s_{i}, t_{j}}\right)\\
&&+  \sum_{i=0}^{m-1} \sum_{j=0} ^{m-1}\left( W_{s_{i}, t_{j+1} }-W_{s_{i}, t_{j}}\right) \left( W_{s_{i+1}, t_{j} }W_{s_{i}, t_{j}}\right).
\end{eqnarray*}

The first summand above usually converges to the stochastic integral with respect to the Wiener sheet $\int_{0}^{s}\int_{0}^{t} W_{u,v} dW_{u, v}$ (see \cite{Imk} or \cite{N81} for the definition of the integral and the proof of the convergence; but the stochastic integral can be also understood as the Skorohod integral $\delta \left(W _{\cdot, \cdot }1_{[0,s]\times [0,t] }(\cdot , \cdot)\right)$ with respect to the Gaussian process $W$). The summand $\sum_{i=0}^{m-1} \sum_{j=0} ^{m-1}\left( W_{s_{i}, t_{j+1} }-W_{s_{i}, t_{j}}\right) \left( W_{s_{i+1}, t_{j} }W_{s_{i}, t_{j}}\right)$ is a specific two-parameter case term and it converges, for fixed $s,t \in [0,1]$  as $m,n\to \infty$ to a random variable denoted $\tilde{M}_{s,t}$ such that the stochastic process $(\tilde{M}_{s,t})_{s,t\in [0,1]} $ is a two-parameter martingale (we refer again to \cite{Imk} or \cite{N81} for details). A similar term $\tilde{M}$ appears also in the stochastic calculus for the fractional Brownian sheet (see \cite{TV}). The stochastic process $\tilde{M}$ can be interpreted as
$$\tilde{M}_{s,t}= \int_{0}^{s} \int_{0} ^{t} d_{u} W_{u,v} d_{v} W_{u,v}.$$
We refer to \cite{TV} for a discussion in this sense.
In terms of the divergence integral it can be written as
$$\tilde {M} _{s,t} = \delta _{(u_{1}, v_{1}), (u_{2}, v_{2})}^{(2)} \left(1_{[0,s] \times [0,t] } ^{\otimes 2} ((u_{1}, v_{1}), (u_{2}, v_{2}) ) 1_{[0,u_{1} ]}(u_{2} ) 1_{[0, v_{2} ]} (v_{1}))\right).$$
Let  us explain the above notation. The symbol $\delta ^{(2)} $ means the double (iterated) Skorohod integral with respect to the Wiener sheet $W$. The indices $(u_{1}, v_{1}), (u_{2}, v_{2})$ in the  formula $\delta _{(u_{1}, v_{1}), (u_{2}, v_{2})}^{(2)}$ means that the double Skorohod integral acts with respect to the two variables $(u_{1}, v_{1})$ and $(u_{2}, v_{2})$.

 More generally, if $(X_{s,t})_{s,t\in [0,1]}  $ is an adapted two-parameter process (let us be vague on this concept because several definitions of the adaptability exist in the two-parameter context; we refer again to \cite{Imk}, \cite{N81} or \cite{Dozzi}), then the sequence
\begin{equation*}
\sum_{i=0}^{m-1} \sum_{j=0}^{n-1} X_{s_{i}, t_{j}}\left( W_{s_{i}, t_{j+1} }-  W_{s_{i}, t_{j} }\right) \left( W_{s_{i+1}, t_{j} }-  W_{s_{i}, t_{j} }\right)
\end{equation*}
converges to the integral
\begin{equation*}
\int_{0}^{t} \int_{0}^{s}  X_{u,v} dM_{u,v} = \int_{0}^{t}\int_{0} ^{s} X_{u,v} d_{u}W_{u,v} d_{v} W_{u,v}
\end{equation*}
which can be also written as a double Skokorod integral in the following way (assuming that $X$ is regular enough in the Malliavin calculus sense)
\begin{eqnarray*}
\int_{0}^{t}\int_{0} ^{s} X_{u,v} d_{u}W_{u,v} d_{v} W_{u,v}&=& \delta ^{(2)} \left( X_{u_{1},v_{2}} 1_{[0,u_{1}]}(u_{2}) 1_{[0,v_{2}] }(v_{1})\right)
\end{eqnarray*}
where $\delta ^{(2)} $ denotes the double Skorohod integrals which as above acts  with respect to the variables $(u_{1}, v_{1}), (u_{2}, v_{2}) \in [0,1] ^{2}$.

\vskip0.2cm

The above discussion motivate the following definition  of the 2D stochastic current.

\begin{definition}\label{defxi}
Let $(W_{s,t}) _{s,t\in [0,1]} $ be a Wiener sheet and for every $x\in \mathbb{R}$, $s,t\in [0,1]$ consider the functional $\Delta (x-W_{s,t}) $ which is a distribution in the Watanabe sense. We define the stochastic current over  the Wiener sheet $W$ as the mapping
\begin{equation}\label{xi}
x\to \xi (x):= \delta ^{(2)} _{(u_{1}, v_{1}), (u_{2}, v_{2}) } \left( \Delta (x-W_{u_{1}, v_{2}}) 1_{[0, u_{1}] }(u_{2}) 1_{[0, v_{2}]} (v_{1})\right).
\end{equation}
\end{definition}

The meaning of this definition will be clarified in the sequel. The term $\Delta (x-W_{u_{1}, v_{2}})$ is a distribution in the Sobolev-Watanabe spaces. In the rest of our paper we will analyze the mapping $\xi$ introduced above from different perspectives: first as a functional of the randomness
 $\omega$ and then as a function of $x\in \mathbb{R} ^{d}$.



\section{The 2D-stochastic current in Watanabe spaces}
Let us first analyze the function $\xi (x)$, for fixed $x\in \mathbb{R}$, as a distribution in the Watanabe (or Sobolev-Watanabe) sense (see Section \ref{prel}). Concretely, we will prove in this section that for every $x$ the function $\xi (x)$ belongs to the Watanabe space $\mathbb{D} ^{\alpha, 2} $ for every $\alpha < -\frac{1}{2}$. We also prove that $\xi (x)$ can be approximated, under the norm of the space $\mathbb{D} ^{\alpha, 2}$ by some Riemann sums. In this way, we make the link with the discussion in the previous section.

\subsection{Regularity of the 2D stochastic current in the Watanabe spaces}


 Our method to prove the regularity of the functional (\ref{xi}) is based on the Wiener-It\^o chaos
decomposition. Let $(X_{s})_{s\in S}$ be an isonormal Gaussian process with covariance
$$EX_{s}X_{t} :=R(s,t)\hskip0.5cm  s,t \in S.$$
The set $S$ is a subset of $\mathbb{R} ^{N}$. We will use the following decomposition of the delta Dirac
function (see Nualart and Vives \cite{NV}, Imkeller et al. \cite{Imk1}, Eddahbi et al. \cite{Ed}) into orthogonal sum of  multiple Wiener-It\^o
integrals
\begin{eqnarray}
\label{deltagen} \Delta (x-X_{s}) &=& \sum _{n\geq 0} R(s)
^{-\frac{n}{2}}p_{R(s)}(x) H_{n} \left(
\frac{x}{R(s)^{\frac{1}{2}}}\right) I_{n}\left( 1_{[0,s]}^{\otimes
n}\right)\\
&: =& \sum _{n\geq 0} a_{n}^{x}(s) I_{n}\left( 1_{[0,s]}^{\otimes
n}\right), \hskip0.5cm x\in \mathbb{R}, s\in S \nonumber
\end{eqnarray}
where $R(s):=R(s,s)$, $p_{R(s)}$ is the Gaussian kernel of variance
$R(s)$, $H_{n}$ is the Hermite polynomial of degree $n$ and $I_{n}$
represents the multiple Wiener-It\^o integral of degree $n$ with
respect to the Gaussian process $X$ as defined in the previous section. We also denoted by
\begin{equation}
\label{a1}
a_{n}^{x}(s)=  R(s)
^{-\frac{n}{2}}p_{R(s)}(x) H_{n} \left(
\frac{x}{R(s)^{\frac{1}{2}}}\right)
\end{equation}
for every $s\in S, x\in \mathbb{R}, n\geq 0$.

For $X=W= (W_{s,t}) _{s,t\in [0,1]}$, the Wiener sheet, we have
\begin{equation}\label{deltasheet}
\Delta (x-W_{u_{1}, v_{2}}) = \sum_{m\geq 0} a_{m}^{x} (u_{1}, v_{2}) I_{m} \left(   1_{[0,u_{1}]\times [0, v_{2}]}^{\otimes m} \right)
\end{equation}
where for every $u_{1}, v_{2} \in [0,1]$,
\begin{equation}\label{amx}
a_{m}^{x} (u_{1}, v_{2}):=(u_{1}v_{2}) ^{-\frac{m}{2}} p_{u_{1}v_{2}} (x) H_{m} \left( \frac{x}{\sqrt{u_{1}v_{2} } } \right) .
\end{equation}

We will show in the next result that the 2D-stochastic current defined in Definition \ref{defxi} is well-defined as a distribution in the Watanabe sense.

\begin{theorem}\label{th1}
For every $x\in \mathbb{R}$ the functional $\xi (x)$ defined in Definition \ref{defxi} belongs to the Watanabe space $\mathbb{D} ^{-\alpha ,2} $ for any $\alpha > \frac{1}{2}$.

\end{theorem}
{\bf Proof: } Let $x\in \mathbb{R}$ be fixed throughout the proof. We will compute first the divergence integral
$$\delta ^{(2)} _{(u_{1}, v_{1}), (u_{2}, v_{2}) } \left( \Delta (x-W_{u_{1}, v_{2}}) 1_{[0, u_{1}] }(u_{2}) 1_{[0, v_{2}]} (v_{1})\right)$$
with respect to the variables $(u_{1}, v_{1})$ and $(u_{2}, v_{2})$ belonging to $[0,1] ^{2}$. Using (\ref{deltasheet}) and the chaotic form of the divergence integral (\ref{sko}), we will  need to symmetrize the function
\begin{eqnarray*}
&&\left( (x_{1}, y_{1}), (x_{2}, y_{2}),\ldots , (x_{m+1}, y_{m+1}), (x_{m+2}, y_{m+2}) \right)\\
 &\to& a_{m} ^{x} (x_{m+1},y_{m+2}) 1_{[x_{m+2}, 1] } (x_{m+1} )1_{[0,y_{m+2}] } (y_{m+1} )  1_{[0, x_{m+1} ] \times [0, y_{m+2}]} ^{\otimes m}\left( (x_{1}, y_{1}), (x_{2}, y_{2}),\ldots , (x_{m}, y_{m}) \right)
\end{eqnarray*}
with respect to the $m+2$ variables
\begin{equation*}
(x_{1}, y_{1}), (x_{2}, y_{2}),\ldots , (x_{m+}, y_{m+2}).
\end{equation*}
The above function is already symmetric in its first $m$ variables. Noting that for any function $f(a_{1}, a_{2},\ldots , a_{m}, x,y)$ which is symmetric with respect to $a_{1}, .., a_{m}$, its symmetrization with respect to its all $m+2$ variables is given by
\begin{eqnarray}\label{sim2}
\tilde{f}(a_{1},.., a_{m+2})= \frac{1}{(m+1)(m+2)} \sum_{k,l=1; k\not=l}^{m+2}f(a_{1},..,\hat{a_{k}}, ..,\hat{a_{l}}, ..a_{m+2}, a_{k}, a_{l})
\end{eqnarray}
(the notation $\hat{a}$ means that the variable $a$ is missing above)  we get, for every $x\in \mathbb{R}$,
\begin{eqnarray}
\xi(x)
&=&
\sum_{m\geq 0}\frac{1}{(m+1)(m+2)}I_{m+2}\sum_{k,l=1; k\not= l} ^{m+2} a_{m} ^{x} (x_{k}, y_{l}) 1_{[0, x_{k}]\times [0, y_{l} ] } (x_{l}, y_{k})\nonumber \\
  &&\times 1_{[0, x_{k} ] \times [0, y_{l}]} ^{\otimes m}\left( (x_{1}, y_{1}), (x_{2}, y_{2}),\ldots , \overline{ (x_{k}, y_{k}) }, \overline{ (x_{l}, y_{l}) },  \ldots ,(x_{m+2}, y_{m+2}) \right) \nonumber \\
  &=& \sum_{m\geq 0} I_{m+2}\frac{1}{(m+1)(m+2)}\sum_{k,l=1; k\not= l} ^{m+2} a_{m} ^{x} (x_{k}, y_{l}) \nonumber\\
  &&1_{[0,x_{k}]}^{\otimes m+1} (x_{1}, .., \overline{x_{k}},.., x_{m+2} )1_{[0,y_{l}] }^{\otimes m+1} (y_{1},.., \overline{y_{l}},.., y_{m+2} )\label{xi-sim}
\end{eqnarray}
where $\overline{ (x_{k}, y_{k}) }$ means, as above,  that the variable $(x_{k}, y_{k})$ is missing.

\b Using (\ref{d12}),  the $\mathbb{D} ^{\alpha ,2}$ norm of the random distribution $\xi (x)$ can be written as
\begin{eqnarray*}
&&\Vert \xi (x) \Vert _{\alpha , 2} ^{2} \\
&=& \sum_{m\geq 0} (1+m) ^{\alpha } (m+2)! \Vert   \left( \frac{1}{m+1}\frac{1}{m+2} \sum_{k,l=1; k\not= l} ^{m+2} a_{m} ^{x} (x_{k}, y_{l}) 1_{[0, x_{k}]\times [0, y_{l} ] } (x_{l}, y_{k}) \right.\\
  &&\left. \times 1_{[0, x_{k} ] \times [0, y_{l}]} ^{\otimes m}\left( (x_{1}, y_{1}), (x_{2}, y_{2}),\ldots , \overline{ (x_{k}, y_{k}) }, \overline{ (x_{l}, y_{l}) },  \ldots ,(x_{m+2}, y_{m+2}) \right)\right) \Vert _{L^{2} (([0,1] ^{2}) ^{m+2} )}\\
&=& \sum_{m\geq 0} (1+m) ^{\alpha } (m+2)! \left( \frac{1}{m+1}\frac{1}{m+2}\right) ^{2} \sum_{k, l=1; k\not= l}^{m+2} \sum_{k', l'=1; k'\not= l'}^{m+2} \int_{([0, 1] ^{2} ) ^{m+2} } dx_{1}dy_{1}..dx_{m+2} dy_{m+2} \\
&&\sum_{k,l=1; k\not= l} ^{m+2} a_{m} ^{x} (x_{k}, y_{l})1_{[0,x_{k}]}^{\otimes m+1} (x_{1}, .., \overline{x_{k}},.., x_{m+2} )1_{[0,y_{l}] }^{\otimes m+1} (y_{1},.., \overline{y_{l}},.., y_{m+2} )  \\
&&\sum_{k',l'=0; k'\not= l'} ^{m+2} a_{m} ^{x} (x_{k'}, y_{l'}) 1_{[0,x_{k'}]}^{\otimes m+1} (x_{1}, .., \overline{x_{k'}},.., x_{m+2} )1_{[0,y_{l'}] }^{\otimes m+1} (y_{1},.., \overline{y_{l'}},.., y_{m+2} ).
\end{eqnarray*}
Notice that if $k\not= k'$,   the terms  $1_{[0, x_{k} ] } (x_{k'}) $ and $1_{[0, x_{k'} ] } (x_{k}) $ appear in the above expression. Therefore, all the summand with $k\not= k'$ will vanish. Similarly, the summands with $l\not= l'$ will vanish. Consequently,
\begin{eqnarray*}
\Vert \xi (x) \Vert _{\alpha , 2} ^{2} &=&\sum_{m\geq 0} (1+m) ^{\alpha } m!\frac{1}{m+1}\frac{1}{m+2}  \sum_{k,l=1; k\not= l}^{m+2} \int_{([0, 1] ^{2} ) ^{m+2} } dx_{1}dy_{1}..dx_{m+2} dy_{m+2} \\
&&\times (a_{m}^{x} (x_{k}, y_{l}))^{2}  1_{[0,x_{k}]}^{\otimes m+1} (x_{1}, .., \overline{x_{k}},.., x_{m+2} )1_{[0,y_{l}] }^{\otimes m+1} (y_{1},.., \overline{y_{l}},.., y_{m+2} )\\
&=&\sum_{m\geq 0} (1+m) ^{\alpha } m!\frac{1}{m+1}\frac{1}{m+2}  \sum_{k,l=1; k\not= l} ^{m+2}\int_{0}^{1} dx_{k} \int_{0}^{1}dy_{l}(a_{m}^{x} (x_{k}, y_{l}))^{2}  \\
&&\times \int_{[0,x_{k}]^{m+1}} dx_{1}...\overline{dx_{k}} ..dx_{m+2} \int_{[0,y_{l}]^{m+1}} dy_{1}...\overline{dy_{l}} ..dy_{m+2}\\
&=&\sum_{m\geq 0} (1+m) ^{\alpha } m!\frac{1}{m+1}\frac{1}{m+2}  \sum_{k,l=1; k\not= l}^{m+2} \int_{0}^{1} dx_{k} \int_{0}^{1}dy_{l}(a_{m}^{x} (x_{k}, y_{l}))^{2}x_{k}^{m+1}y_{l} ^{m+1}.
\end{eqnarray*}
Using the expression (\ref{amx}) of the coefficient $a_{m}^{x}$ we have
\begin{equation*}
(a_{m}^{x} (x_{k}, y_{l}))^{2}x_{k}^{m+1}y_{l} ^{m+1} = \left( p_{x_{k}y_{l}}(x) H_{m} \left( \frac{x}{\sqrt{x_{k}y_{l}}} \right)\right) ^{2} x_{k}y_{l}.
\end{equation*}
Therefore
\begin{eqnarray*}
\Vert \xi (x) \Vert _{\alpha , 2} ^{2}=\sum_{m\geq 0} (1+m) ^{\alpha } m!\int_{0}^{1} \int_{0}^{1} dudv \left( p_{uv}(x) H_{m} \left( \frac{x}{\sqrt{uv} } \right) \right)^{2}uv.
\end{eqnarray*}
We use  the identity, for every $y\in \mathbb{R}$,
\begin{equation}
\label{id} H_{n}(y)e^{-\frac{y^{2}}{2}} = (-1)^{[n/2]}
2^{\frac{n}{2}} \frac{2}{n!\pi } \int _{0}^{\infty } u^{n}
e^{-u^{2}} g(uy\sqrt{2}) du
\end{equation}
where $g(r)=\cos (r)$ if $n$ is even and $g(r)=\sin (r)$ if $n$ is
odd. Since $\vert g(r)\vert \leq 1$, we have the bound
\begin{equation}\label{cn}
\left|H_{n}(y)e^{-\frac{y^{2}}{2}} \right| \leq
2^{\frac{n}{2}}\frac{2}{n!\pi } \Gamma (\frac{n+1}{2}):=c_{n}.
\end{equation}
This implies that
\begin{eqnarray*}
\Vert \xi (x) \Vert _{\alpha , 2} ^{2}&\leq &\frac{1}{2\pi}\sum_{m\geq 0} (1+m) ^{\alpha } m!c_{m}^{2} \int_{0}^{1} \int_{0}^{1} dudv (uv) ^{\frac{1}{2}}\\
&=& \frac{1}{2\pi}\frac{4}{9}\sum_{m\geq 0} (1+m) ^{\alpha } m!c_{m}^{2}
\end{eqnarray*}
and this is finite for every $\alpha <\frac{-1}{2}$ because, by Stirling's formula,  $m!c_{m}^{2}$ (this has been already used in \cite{NV}, \cite{Ed}, \cite{CNT}) behaves as $n^{-\frac{1}{2}}$. \qed

\begin{remark}
We remark that the regularity of the functional $\xi (x)$ in the Watanabe sense is the same as in the one-parameter case (see \cite{FT}).
This implies in particular that it is not a random variable, but a distribution.
\end{remark}

\begin{remark}
Note that the integration interval for the stochastic integral $dW$ in the definition of $\xi(x)$ is the entire interval $[0,1]$. In the one-parameter case (see \cite{FT}), we need to consider the quantity $\xi(x)= \int_{a} ^{1} \Delta (x-B_{s})dB_{s}$ with $a>0$ in order to obtain that $\xi (x) \in \mathbb{D}^{\alpha , 2}$ with $\alpha <\frac{-1}{2}$.
\end{remark}

\subsection{Approximation of 2D-stochastic current in Watanabe spaces}
\vskip0.3cm
Let us introduce the Riemann sum
\begin{equation}
\label{TNM}
T_{N,M} (x)= \sum_{i=0} ^{N-1} \sum_{j=0} ^{M-1} \Delta (x-W_{s_{i}, t_{j}} ) \left( W_{s_{i+1}, t_{j}}-W_{s_{i}, t_{j}}\right) \left( W_{s_{i}, t_{j+1}}-W_{s_{i}, t_{j}}\right).
\end{equation}

\begin{remark}
 Following the lines of Theorem \ref{th1}, we can prove that $T_{N,M}(x)$ belongs to the Sobolev-Watanabe space $\mathbb{D} ^{\alpha, 2 }$ for every $\alpha <\frac{-1}{2}$. This is also a consequence of the proof of Theorem \ref{tt2} below.
\end{remark}

Following the discussion at the beginning of Section \ref{defint} it is natural to expect that $T_{N,M}(x)$  converges as $M,N\to \infty$ to $\xi (x)$ for every $x\in \mathbb{R}$. We will actually show that the sequence (\ref{TNM}) approximates $\xi (x)$ in the norm of the Watanabe space $\mathbb{D} ^{-\alpha ,2} $ with $\alpha >\frac{1}{2}$.

\begin{theorem}\label{tt2}
For every $x\in \mathbb{R}$, the sequence $T_{N,M} (x)$ given by (\ref{TNM}) converges to $\xi (x)$ as $M, N\to \infty$ in the norm $\Vert \cdot \Vert _{-\alpha ,2}$ for any $\alpha >\frac{1}{2}$.
\end{theorem}
{\bf Proof: } Fix again $x\in \mathbb{R}$. The first step is to express $T_{N,M}$  using multiple Wiener-It\^o integrals. We have, by the integration by parts formula (\ref{ip}) and  the symmetrization formula (\ref{sim2})
\begin{eqnarray}
T_{N,M} (x)
&=& \delta ^{(2)} \left(\sum_{i=0}^{N-1} \sum_{j=0}^{M-1} \Delta (x-W_{s_{i},t_{j}})1_{[0,s_{i}]\times [t_{j}, t_{j+1}]} (u_{1},v_{1}) 1_{[s_{i}, s_{i+1}]\times [0, t_{j}]} (u_{2}, v_{2})\right)\nonumber \\
&=& \sum_{i=0}^{N-1}\sum_{j=0}^{M-1}  a_{m}^{x}(s_{i},t_{j}) \nonumber\\
&&\times \sum_{m\geq 0} I_{m+2} \left( 1_{[0,s_{i} ]\times [0,t_{j}]}^{\otimes m} ((x_{1},y_{1}),..., (x_{m},y_{m}))  1_{[0,s_{i}]\times [t_{j}, t_{j+1}]} (u_{1},v_{1}) 1_{[s_{i}, s_{i+1}]\times [0, t_{j}]} (u_{2}, v_{2}) \right)\nonumber\\
&=&\sum_{i=0}^{N-1}\sum_{j=0}^{M-1}  a_{m}^{x}(s_{i},t_{j})\sum_{m\geq 0} I_{m+2}\left( \frac{1}{(m+1)(m+2)} \right. \nonumber\\
&& \left. \sum_{k,l=1; k\not= l}^{m+2} 1_{[0, s_{i}]\times [0, t_{j}]} ^{\otimes m}((x_{1},y_{1}),.., \overline{(x_{k},y_{k})}, \overline{(x_{l},y_{l})},.., (x_{m+2},y_{m+2}))\right. \nonumber\\
&&\left. 1_{[0,s_{i}]\times [t_{j}, t_{j+1}]} (x_{k}, y_{k}))1_{[s_{i}, s_{i+1}]\times [0, t_{j}]} (x_{l}, y_{l})\right).\label{tnmsim}
\end{eqnarray}
Taking into account relation (\ref{tnmsim}) and its analogous for $\xi (x)$ (see the proof of Theorem \ref{th1}, formula \ref{xi-sim})  we obtain

\begin{eqnarray*}
\Vert T_{N,M}(x)-\xi (x) \Vert ^{2} _{\alpha , 2}
 &=&\sum_{m\geq 0} (3+m) ^{\alpha } \frac{m!} {(m+1) (m+2)} \sum_{k,l=1; k\not= l}^{m+2}
\int_{[0,1] ^{2m}} dx_{1}..dx_{m+2} dy_{1}..dy_{m+2} \\
&& \left(  \sum_{i=0}^{N-1} \sum_{j=0}^{M-1}a_{m}^{x}(s_{i}, t_{j}) 1_{[0,s_{i}]} ^{\otimes m+1} (x_{1},..,\overline{x_{k}},.., x_{m+2})1_{[s_{i},s_{i+1}]}(x_{k}) 1_{[0,t_{j}]} ^{\otimes m+1} (y_{1},..,\overline{y_{l}} ,..y_{m+2}) \right. \\
&& \left. -  a_{m}^{x} (x_{k}, y_{l}) 1_{[0,x_{k} ]}^{\otimes m+1} (x_{1},..,\overline{x_{k}},.., x_{m+2})1_{[0,y_{l}]}^{\otimes m+1} (y_{1},..,\overline{y_{l}} ,..y_{m+2}) \right)^{2}
\end{eqnarray*}
and as in the proof of Theorem \ref{th1}
\begin{eqnarray*}
\Vert T_{N,M}(x)-\xi (x) \Vert ^{2} _{\alpha , 2}
&=&\sum_{m\geq 0} (3+m) ^{\alpha }\frac{m!} {(m+1) (m+2)} \sum_{k,l=1; k\not= l}^{m+2} \\
 &&\left( \sum_{i=0}^{N-1} \sum_{j=0}^{M-1}(a_{m}^{x}(s_{i}, t_{j}))^{2} s_{i}^{m+1} (s_{i+1}-s_{i}) t_{j}^{m+1} (t_{j+1}-t_{j})  \right. \\
&&\left. - 2 \sum_{i=0}^{N-1} \sum_{j=0}^{M-1}a_{m}^{x}(s_{i}, t_{j})s_{i}^{m+1}  t_{j}^{m+1}\int_{s_{i}}^{s_{i+1}} \int_{t_{j}}^{t_{j+1}}dx_{k}dy_{l}a_{m}^{x} (x_{k}, y_{l})  \right. \\
&&\left. + \int_{0}^{1} \int_{0}^{1} dx_{k}dy_{l} (a_{m}^{x}(x_{k}, y_{l})) ^{2} x_{k}^{m+1} y_{l} ^{m+1} \right)\\
&=&\sum_{m\geq 0} (3+m) ^{\alpha }m! \left( \sum_{i=0}^{N-1} \sum_{j=0}^{M-1}(a_{m}^{x}(s_{i}, t_{j}))^{2} s_{i}^{m+1} (s_{i+1}-s_{i}) t_{j}^{m+1} (t_{j+1}-t_{j})  \right. \\
&&\left. -  2 \sum_{i=0}^{N-1} \sum_{j=0}^{M-1}a_{m}^{x}(s_{i}, t_{j})s_{i}^{m+1}  t_{j}^{m+1}\int_{s_{i}}^{s_{i+1}} \int_{t_{j}}^{t_{j+1}}dudv a_{m}^{x} (u,v)\right. \\
&&\left.  +  \int_{0}^{1} \int_{0}^{1}dudv (a_{m}^{x}(u,v))^{2} (uv) ^{m+1}\right).
\end{eqnarray*}
Now we claim that, for every fixed $m\geq 0$, the sequence
\begin{eqnarray*}
&&\left( \sum_{i=0}^{N-1} \sum_{j=0}^{M-1}(a_{m}^{x}(s_{i}, t_{j}))^{2} s_{i}^{m+1} (s_{i+1}-s_{i}) t_{j}^{m+1} (t_{j+1}-t_{j})  \right. \\
&&\left. -  2 \sum_{i=0}^{N-1} \sum_{j=0}^{M-1}a_{m}^{x}(s_{i}, t_{j})s_{i}^{m+1}  t_{j}^{m+1}\int_{s_{i}}^{s_{i+1}} \int_{t_{j}}^{t_{j+1}}dudv a_{m}^{x} (u,v) + 2 \int_{0}^{1} \int_{0}^{1}dudv (a_{m}^{x}(u,v))^{2} (uv) ^{m+1}\right)
\end{eqnarray*}
converges to zero as $N,M\to \infty$. This is true because
\begin{equation*}
\sum_{i=0}^{N-1} \sum_{j=0}^{M-1}(a_{m}^{x}(s_{i}, t_{j}))^{2} s_{i}^{m+1} (s_{i+1}-s_{i}) t_{j}^{m+1} (t_{j+1}-t_{j}) \to _{N,M\to \infty} \int_{0}^{1} \int_{0}^{1}dudv (a_{m}^{x}(u,v))^{2} (uv) ^{m+1}
\end{equation*}
and
\begin{equation*}
\sum_{i=0}^{N-1} \sum_{j=0}^{M-1}a_{m}^{x}(s_{i}, t_{j})s_{i}^{m+1}  t_{j}^{m+1}\int_{s_{i}}^{s_{i+1}} \int_{t_{j}}^{t_{j+1}}dudv a_{m}^{x}(u,v)\to _{N,M\to \infty} \int_{0}^{1} \int_{0}^{1}dudv (a_{m}^{x}(u,v))^{2} (uv) ^{m+1}
\end{equation*}
using Riemann sums convergence.
To conclude the convergence of $\Vert T_{N,M}(x)-\xi (x) \Vert ^{2} _{\alpha , 2}$ to zero for every $\alpha <\frac{-1}{2}$ it suffices to check that $$\sum_{i=0}^{N-1} \sum_{j=0}^{M-1}(a_{m}^{x}(s_{i}, t_{j}))^{2} s_{i}^{m+1} (s_{i+1}-s_{i}) t_{j}^{m+1} (t_{j+1}-t_{j})\leq cm^{-\frac{1}{2}},$$

$$\sum_{i=0}^{N-1} \sum_{j=0}^{M-1}a_{m}^{x}(s_{i}, t_{j})s_{i}^{m+1}  t_{j}^{m+1}\int_{s_{i}}^{s_{i+1}} \int_{t_{j}}^{t_{j+1}}dudv a_{m}^{x} (u,v) \leq cm^{-\frac{1}{2}},$$

and
$$ \int_{0}^{1} \int_{0}^{1}dudv (a_{m}^{x}(u,v))^{2} (uv) ^{m+1}\leq cm^{-\frac{1}{2}}.$$
with $c>0$ a constant not depending on $m$. The last bound has been proved before (see the proof of Theorem \ref{th1}) and the first two bounds can be obtained analogously using relations (\ref{id}) and (\ref{cn}).\qed

\section{2D-stochastic currents in deterministic Sobolev spaces}

The purpose is this paragraph in to study the mapping (\ref{xi}) as a function of the spatial variable $x\in \mathbb{R}$.

\subsection{Regularity in the deterministic Sobolev spaces}
We study in this part  the regularity with respect to
the  variable $x  \in \mathbb{R}^{d} $ of the mapping  given by
(\ref{xi})  in the (deterministic) Sobolev spaces $H^{-r}
(\mathbb{R}^{d}; \mathbb{R}^{d})$.
This is the dual space to $H^{r}\left(  \mathbb{R}^{d},\mathbb{R}^{d}\right)
$, or equivalently the space of all vector valued distributions $\varphi$ such that
$\int_{\mathbb{R}^{d}}\left(  1+\left|  x\right|  ^{2}\right)  ^{-r}\left|
\widehat{\varphi}\left(  x\right)  \right|  ^{2}dx<\infty$, $\widehat{\varphi
}$ being the Fourier transform of $\varphi$.
Assume first that $d=1$. The case $d\geq 2$ will be treated later.

We will recall the following lemma (see \cite{FT}, Lemma 1).

\begin{lemma}\label{ahat}
Let $s=(s_{1},.., s_{N})\in [0,1]^{N}$ be fixed. The Fourier transform of the function
\begin{equation*}
x\to a_{m} ^{x} (s) : \mathbb{R} \to \mathbb{R}
\end{equation*}
is denoted by $a_{m}^{\hat{x}}(s)$ and it is given by
\begin{equation*}
a_{m}^{\hat{x}}(s)= e^{-\frac{x^{2}}{2} \vert s\vert } \frac{ (-i) ^{m} x^{m} }{m!}
\end{equation*}
where $\vert s\vert =s_{1}...s_{N}$.
\end{lemma}

As a consequence, the Fourier transform with respect to $x$, denoted in the sequel $\Delta (\widehat{x}-W_{u_{1}, v_{2}}) $, of the random distribution $x\in \mathbb{R} \to \Delta (x-W_{u_{1}, v_{2}}) $ is given by
\begin{equation*}
\Delta (\widehat{x}-W_{u_{1}, v_{2}})=\sum_{m} a_{m}^{\hat{x}}(u_{1}, v_{2}) I_{m} \left( 1_{[0,u_{1}]\times [0, v_{2}]} ^{\otimes m} \right)= e^{-ixW_{u_{1}, v_{2}} }.
\end{equation*}

\begin{prop}\label{pp1}
The Fourier transform of the function
\begin{equation*}
x\to \xi (x)
\end{equation*}
is given by
\begin{equation*}
\hat{\xi} (x) = \delta _{(u_{1},v_{1}), (u_{2},v_{2})}^{(2)} \left(   e^{-ixW_{u_{1}, v_{2}}}1_{[0, u_{1}] \times [0, v_{2}] } \right).
\end{equation*}
\end{prop}
{\bf Proof: }Using formula (\ref{xi-sim}) and the expression of the Fourier transform of $x\to a_{m}^{x}$ from Lemma \ref{ahat}, we get (to justify the interchange of the order of integration below we refer to Exercise 3.2.7 in \cite{N})
\begin{eqnarray*}
\hat{\xi} (x) &= &\sum_{m\geq 0} I_{m+2}\frac{1}{(m+1)(m+2)}\sum_{k,l=1; k\not= l} ^{m+2} a_{m} ^{\widehat{x}} (x_{k}, y_{l}) \nonumber\\
  &&1_{[0,x_{k}]}^{\otimes m+1} (x_{1}, .., \overline{x_{k}},.., x_{m+2} )1_{[0,y_{l}] }^{\otimes m+1} (y_{1},.., \overline{y_{k}},.., y_{m+2} )\\
  &=& \delta ^{(2)} _{(u_{1}, v_{1}), (u_{2}, v_{2})}  e^{-\frac{x^{2}}{2}u_{1}v_{2}}\sum_{m\geq 0} \frac{(-i)^{m} x^{m}}{m!} I_{m}\left(1^{\otimes m}_{[0, u_{1}] \times [0,v_{2}] } (\cdot ) \right) \\
  &=& \delta ^{(2)} _{(u_{1}, v_{1}), (u_{2}, v_{2})} \left( e^{-ixW_{u_{1}, v_{2}}}1_{[0, u_{1}] \times [0, v_{2}] } \right).
\end{eqnarray*}
Here we used the fact that
\begin{equation*}
e^{-ixW_{u_{1}, v_{2}} }= \sum_{m\geq 0} \frac{(-i)^{m} x^{m}}{m!} I_{m}\left(1^{\otimes m}_{[0, u_{1}] \times [0,v_{2}] } (\cdot )\right)
\end{equation*}
which can be obtained using Stroock's formula (see \cite{FT} relation (9)). \qed
\vskip0.3cm

Consequently, we obtain the following result:
\begin{theorem}\label{th3}
For almost all $\omega$, it holds that $\xi \in H^{-r} (\mathbb{R}, \mathbb{R}) $  if  $r>\frac{1}{2}$
\end{theorem}
{\bf Proof: }It is easy to see that
\begin{eqnarray*}
\mathbf{E}\left| \widehat{\xi }(x) \right| ^{2}&=&\int_{0} ^{1} du_{1} \int_{0} ^{u_{1} } du_{2} \int_{0}^{1}dv_{2} \int_{0} ^{v_{2}} dv_{1} \left| e^{-ixW_{u_{1}, v_{2}} }\right| ^{2}\\
&=& \int_{0}^{1} u_{1} du_{1} \int_{0} ^{1} v_{2}dv_{2}=\frac{1}{4}
\end{eqnarray*}
and consequently
\begin{eqnarray*}
\mathbf{E} \Vert \xi(x) \Vert _{H^{-r} }=\int_{\mathbb{R}} \frac{1}{(1+ x^{2})^{r}} \mathbf{E} \left| \widehat{\xi }(x) \right| ^{2}dx
\end{eqnarray*}
is finite if and only if $r>\frac{1}{2}$. \qed

\begin{remark}
We note that, as in the case of 1D -currents, the regularity of (\ref{xi}) is the same with respect to the space variable $x$ and with respect to the randomness variable $\omega$.
\end{remark}
\subsection{Convergence of the Fourier transforms}

Here we study the convergence of the sequence $T_{N,M}$ (\ref{TNM}) to the 2D-stochastic current (\ref{xi}) in the Sobolev space. Precisely, we have the following:

\begin{theorem}
Let $T_{N,M}$ be given by (\ref{TNM}). Then as $M,N\to \infty$, the sequence $T_{N,M}$ converges to $\xi$ in $L^{2}(\Omega; H^{-r}(\mathbb{R}; \mathbb{R}))$ for every $r>\frac{1}{2}$.
\end{theorem}
{\bf Proof: } From (\ref{TNM}) it is immediate that
\begin{equation*}
T_{N,M}(\hat{x}) =\sum_{i=0}^{N-1} \sum_{j=0}^{M-1} e^{_ixW_{s_{i}, t_{j}}} \left( W_{s_{i+1}, t_{j}} -  W_{s_{i}, t_{j}}\right) \left(  W_{s_{i}, t_{j+1}}- W_{s_{i}, t_{j}}\right)
\end{equation*}
where  $T_{N,M}(\hat{x})$ denotes the Fourier transform of
$x\to T_{N,M}(x).$
We also know from Proposition \ref{pp1} that
\begin{equation*}
\xi (\hat{x})=\delta ^{(2)}_{(u_{1},v_{1}), (u_{2}, v_{2}) }\left( e ^{-xW_{u_{1}, v_{2}} } 1_{[0,u_{1}]} (u_{2} ) 1_{[0,v_{2} ] } (v_{1} ) \right).
\end{equation*}
Then the $L^{2}(\Omega; H^{-r}(\mathbb{R}; \mathbb{R}))$-norm of the difference $T_{N,M}-\xi$ can be computed as follows
\begin{eqnarray*}
\mathbf{E} \left| T_{N,M} -\xi \right| ^{2} _{H^{-r} (\mathbb{R}; \mathbb{R})} =\mathbf{E}\int_{\mathbb{R} } \left| T_{N,M}(\hat{x}) -\xi(\hat{x} ) \right| ^{2} (1+x^{2}) ^{-r} dx.
\end{eqnarray*}
Let us first compute
\begin{eqnarray*}
&& \mathbf{E} \left| T_{N,M}(\hat{x}) -\xi(\hat{x} ) \right| ^{2}\\
&=&\mathbf{E} \left|  \sum_{i=0}^{N-1} \sum_{j=0}^{M-1} \cos (-xW_{s_{i},t_{j}})  \left( W_{s_{i+1}, t_{j}} -  W_{s_{i}, t_{j}}\right) \left(  W_{s_{i}, t_{j+1}}- W_{s_{i}, t_{j}}\right) \right. \\
&&\left. -\delta ^{(2)}_{(u_{1},v_{1}), (u_{2}, v_{2}) }\left( \cos (-xW_{u_{1}, v_{2}} ) 1_{[0,u_{1}]} (u_{2} ) 1_{[0,v_{2} ] } (v_{1} ) \right)\right| ^{2} \\
&&+  \mathbf{E} \left|  \sum_{i=0}^{N-1} \sum_{j=0}^{M-1} \sin (-xW_{s_{i},t_{j}}  \left( W_{s_{i+1}, t_{j}} -  W_{s_{i}, t_{j}}\right) \left(  W_{s_{i}, t_{j+1}}- W_{s_{i}, t_{j}}\right) \right. \\
&&\left. - \delta ^{(2)}_{(u_{1},v_{1}), (u_{2}, v_{2}) }\left( \sin (-xW_{u_{1}, v_{2}} ) 1_{[0,u_{1}]} (u_{2} ) 1_{[0,v_{2} ] } (v_{1} ) \right)\right| ^{2}
\end{eqnarray*}
and from the independence of the increments of the Brownian sheet
\begin{eqnarray*}
&& \mathbf{E} \left| T_{N,M}(\hat{x}) -\xi(\hat{x} ) \right| ^{2}\\
&=& \mathbf{E} \sum_{i=0}^{N-1}  \sum_{j=0}^{M-1} \cos (-xW_{s_{i},t_{j}}) ^{2} s_{i}t_{j} (s_{i+1}-s_{i}) (t_{j+1}-t_{j}) \\
&&+ \mathbf{E} \int_{0}^{1} du_{1}\int_{0}^{u_{1}}du_{2}\int_{0}^{1} dv_{2} \int_{0}^{v_{2}} dv_{1} \cos (-xW_{u_{1}, v_{2}} ) ^{2} \\
&& -2\mathbf{E}  \sum_{i=0}^{N-1}  \sum_{j=0}^{M-1}\int_{s_{i}}^{s_{i+1} } du_{1} \int_{0}^{s_{i}} du_{2} \int_{t_{j} } ^{t_{j+1}} dv_{2} \int_{0}^{t_{j}}dv_{1} \cos (-xW_{s_{i},t_{j}})\cos (-xW_{u_{1},v_{2}})\\
&+&  \mathbf{E} \sum_{i=0}^{N-1}  \sum_{j=0}^{M-1} \sin (-xW_{s_{i},t_{j}}) ^{2} s_{i}t_{j} (s_{i+1}-s_{i}) (t_{j+1}-t_{j}) \\
&&+ \mathbf{E} \int_{0}^{1} du_{1}\int_{0}^{u_{1}}du_{2}\int_{0}^{1} dv_{2} \int_{0}^{v_{2}} dv_{1} \sin (-xW_{u_{1}, v_{2}} ) ^{2} \\
&& -2\mathbf{E}  \sum_{i=0}^{N-1}  \sum_{j=0}^{M-1}\int_{s_{i}}^{s_{i+1} } du_{1} \int_{0}^{s_{i}} du_{2} \int_{t_{j} } ^{t_{j+1}} dv_{2} \int_{0}^{t_{j}}dv_{1} \sin (-xW_{s_{i},t_{j}})\sin (-xW_{u_{1},v_{2}}).
\end{eqnarray*}
We will obtain, using the well-known trigonometric formulas $\sin^{2}(x)+ \cos ^{2}(x)=1$ and $\sin(x)\sin(y)+ \cos (x)\cos (y)= \cos (x+y)$
\begin{eqnarray*}
 \mathbf{E} \left| T_{N,M}(\hat{x}) -\xi(\hat{x} ) \right| ^{2}
&=& \mathbf{E} \sum_{i=0}^{N-1}  \sum_{j=0}^{M-1}  s_{i}t_{j} (s_{i+1}-s_{i}) (t_{j+1}-t_{j}) +\left(  \int_{0}^{1} udu \right) ^{2}\\
&&-2\mathbf{E}  \sum_{i=0}^{N-1}  \sum_{j=0}^{M-1}s_{i}t_{j} \int_{s_{i}}^{s_{i+1} } du_{1}\int_{t_{j} } ^{t_{j+1}} dv_{2}\cos \left(  x( W_{s_{i}, t_{j}}-W_{u_{1}, v_{2}}) \right).
\end{eqnarray*}
Since $W_{s_{i}, t_{j}}-W_{u_{1}, v_{2}}$ with $s_{i}<u_{1}$ and $t_{j}<v_{2}$ is a centered Gaussian random variable with variance $u_{1}v_{2}-s_{i}t_{j}$ we get
\begin{equation*}
\mathbf{E}\cos \left(  x( W_{s_{i}, t_{j}}-W_{u_{1}, v_{2}}) \right)= e^{-\frac{ x^{2}(u_{1}v_{2}-s_{i}t_{j})}{2}}
\end{equation*}
and therefore
\begin{eqnarray*}
 \mathbf{E} \left| T_{N,M}(\hat{x}) -\xi(\hat{x} ) \right| ^{2}
&=&   \sum_{i=0}^{N-1}  \sum_{j=0}^{M-1}  s_{i}t_{j} (s_{i+1}-s_{i}) (t_{j+1}-t_{j}) +\left(  \int_{0}^{1} udu \right) ^{2}\\
&&-2  \sum_{i=0}^{N-1}  \sum_{j=0}^{M-1}s_{i}t_{j} \int_{s_{i}}^{s_{i+1} } du_{1}\int_{t_{j} } ^{t_{j+1}} dv_{2}e^{-\frac{ x^{2}(u_{1}v_{2}-s_{i}t_{j})}{2}}.
\end{eqnarray*}
We can write
 \begin{eqnarray*}
&&\mathbf{E} \left| T_{N,M} -\xi \right| ^{2} _{H^{-r} (\mathbb{R}; \mathbb{R})} =\mathbf{E}\int_{\mathbb{R} } \left| T_{N,M}(\hat{x}) -\xi(\hat{x} ) \right| ^{2} (1+x^{2}) ^{-r} dx\\
&=&\int_{\mathbb{R} } dx (1+x^{2}) ^{-r} \left(  \sum_{i=0}^{N-1}  \sum_{j=0}^{M-1} s_{i}t_{j} (s_{i+1}-s_{i}) (t_{j+1}-t_{j}) +\left(  \int_{0}^{1} udu \right) ^{2}  \right)\\
&&-2\int_{\mathbb{R} } dx (1+x^{2}) ^{-r} \left(   \sum_{i=0}^{N-1}  \sum_{j=0}^{M-1}s_{i}t_{j} \int_{s_{i}}^{s_{i+1} } du_{1}\int_{t_{j} } ^{t_{j+1}} dv_{2}e^{-\frac{ x^{2}(u_{1}v_{2}-s_{i}t_{j})}{2}}\right)
\end{eqnarray*}

\b By a Riemmann sum convergence, it is clear that
\begin{equation*}
 \sum_{i=0}^{N-1}s_{i}(s_{i+1}-s_{1}) \to _{N\to \infty} \int_{0}^{1}udu, \hskip0.5cm \sum_{j=0}^{M-1}t_{j}(t_{j+1}-t_{j}) \to_{M\to \infty} \int_{0}^{1}vdv
 \end{equation*}
 and for every fixed $x\in \mathbb{R}$
 \begin{equation*}
 \sum_{i=0}^{N-1}  \sum_{j=0}^{M-1}s_{i}t_{j} e^{\frac{ x^{2}s_{i}t_{j}}{2}}\int_{s_{i}}^{s_{i+1} } du_{1}\int_{t_{j} } ^{t_{j+1}} dv_{2}e^{-\frac{ x^{2}u_{1}v_{2}}{2}}\to _{N,M\to \infty} \int_{0}^{1} \int_{0}^{1} uvdudv  e^{\frac{ -x^{2}uv}{2}}e^{\frac{ x^{2}uv}{2}}= \left( \int_{0}^{1}udu\right) ^{2}
 \end{equation*}
 and this implies that for every $x\in \mathbb{R}$  the sequence $\mathbf{E} \left| T_{N,M}(\hat{x}) -\xi(\hat{x} ) \right| ^{2}$ converges to zero as $N,M\to \infty$. The convergence of $\mathbf{E} \left| T_{N,M} -\xi \right| ^{2} _{H^{-r} (\mathbb{R}; \mathbb{R})} $ to zero follows from the dominated convergence theorem since
\begin{eqnarray*} &&(1+x^{2}) ^{-r} \left(  \sum_{i=0}^{N-1}  \sum_{j=0}^{M-1} s_{i}t_{j} \int_{s_{i}}^{s_{i+1} } du_{1}\int_{t_{j} } ^{t_{j+1}} dv_{2} \left( 1-e^{-\frac{ x^{2}(u_{1}v_{2}-s_{i}t_{j})}{2}} \right)\right)\\
&\leq & (1+x^{2}) ^{-r}\left(  \sum_{i=0}^{N-1}  \sum_{j=0}^{M-1} s_{i}t_{j} \int_{s_{i}}^{s_{i+1} } du_{1}\int_{t_{j} } ^{t_{j+1}} dv_{2} \right)\\
&\leq & (1+x^{2}) ^{-r}
\end{eqnarray*}
which is integrable for every $r>\frac{1}{2}$. \qed

\section{Multidimensional 2D-currents}

Newt we study in this parts $d$ -dimensional stochastic currents over a Wiener sheet.
Let $W= \left( W^{(1)}, W^{(2)},..., W^{(d)}\right) $ be a $d$-dimensional Brownian sheet. That means, for every $i=1,..,d$, the component $W^{(i)}= (W^{(i)} _{s,t})_{s,t\in [0,1]}$ is a Wiener sheet and   the components of $W$ are mutually independent.

We will defined the $d$ -dimensional current over the Wiener sheet as follows.

\begin{definition}\label{defxid}The 2D- stochastic current over $W$ in $\mathbb{R} ^{d}$ is given by the mapping
\begin{equation*}
x \in \mathbb{R} ^{d}\to \xi(x)=(\xi _{i,j}) _{i,j=1,..,d}:= \left( \delta ^{i} _{(u_{1}, v_{1}) } \delta ^{j} _{(u_{2}, v_{2})} \left(   \Delta (x-W_{u_{1}, v_{2}}) 1_{[0,u_{1}]}(u_{2}) 1_{[0,v_{2}]} (v_{1})\right)\right)_{i,j=1,..,d}.
\end{equation*}
\end{definition}
Consequently  $\xi(x)$ is, for every $x\in \mathbb{R} ^{d}$, a $d\times d$ matrix whose components are distributions in Sobolev-Watanabe spaces. The
notation $\delta ^{i} $ refers to the divergence integral with respect to the component $W^{(i)}$ for every $i=1,..,d$.

At this point, we will introduce some notation. Let us denote, for every $i,j=1,..,d$,   by
\begin{equation*}
\Delta _{i,j} (x-W_{u,v}) := \prod _{k=1; k\not= i,j} ^{d}\Delta (x_{k} -W^{(k)}_{u,v}) \hskip0.5cm \forall u,v\in [0,1].
\end{equation*}
Therefore, we can formally write
\begin{equation}\label{d1}
\Delta \left( x-W_{u, v}\right) = \Delta _{i,j} (x-W_{u,v}) \Delta (x_{i}-W^{(i)}_{u,v} )  \Delta (x_{j}-W^{(j)}_{u,v} ).
\end{equation}
We have the following result.

\begin{theorem}
The stochastic current $x\to \xi (x)$ introduced in Definition (\ref{defxid}) belongs to the Sobolev space $H^{-r}(\mathbb{R}^{d}; \mathbb{R} ^{d})$ for every $r>\frac{d}{2}$.
\end{theorem}
{\bf Proof: } Let us consider $i\not= j$. Using (\ref{d1}) and the general chaos expansion formula (\ref{deltagen}) for the delta Dirac functional, the component $\xi _{i,j}$ can be decomposed as follows
\begin{eqnarray*}
\xi_{i,j}(x)&=& \delta ^{i} _{(u_{1}, v_{1}) } \delta ^{j} _{(u_{2}, v_{2}) } \left( \Delta _{i,j} (x-W_{u_{1}, v_{2} }) 1_{[0, u_{1} ]} (u_{2}) 1_{[0, v_{2}] }(v_{1})\right) \\
&&\times \sum_{m_{i}\geq 0} a_{m_{i}} ^{x_{i}}(u_{1}, v_{2})I_{m_{i}}^{i} \left( 1_{[0,u_{1}]\times [0, v_{2}] } ^{\otimes m_{i} } \right)
\sum_{m_{j}\geq 0} a_{m_{j}} ^{x_{j}}(u_{1}, v_{2})I_{m_{j}}^{j} \left( 1_{[0,u_{1}]\times [0, v_{2}] } ^{\otimes m_{j} } \right) .
\end{eqnarray*}
Here $I_{m } ^{i}$ denotes the multiple integral of order $m$ with respect to the Wiener sheet $W^{(i)}$. We will get, after symmetrization, for every $i,j=1,..,d$, $i\not=j$

\begin{eqnarray*}
\xi_{i,j}(x)&=& \sum_{m_{i}\geq 0} \sum_{m_{j}\geq 0} \frac{1}{m_{i}+ 1} \frac{1}{m_{j}+1} \\
&&  I_{m_{i}+1}^{i} I_{m_{j}+1} ^{j} \left( \sum_{k=1} ^{m_{i}+1} \sum_{l=1} ^{m_{j}+1} a_{m_{i}}^{x_{i}}(a_{k}, d_{l})a_{m_{j}}^{x_{j}}(a_{k}, d_{l})\Delta _{i,j}(x-W_{a_{k} , d_{l}}) \right.\\
&& \left. 1_{[0,a_{k}] \times [0, d_{l} ] } ^{\otimes m_{i} } \left( (a_{1}, b_{1}), ..., \overline{(a_{k}, b_{k})},..., (a_{m_{i}+1} , b_{m_{i}+1}) \right) 1_{[0,a_{k}]} (c_{l}) 1_{[0,d_{l}]} (b_{k})\right. \\
&& \left. 1_{[0,a_{k}] \times [0, d_{l} ] } ^{\otimes m_{j} }\left( (c_{1}, d_{1}), ..., \overline{(c_{l}, d_{l})},...,(c_{m_{j}+1}, d_{m_{j}+1}) \right) \right) .
\end{eqnarray*}
Let us comment on the above expression. The multiple integral $I_{m_{i}+1} ^{i}$ acts with respect to the variables $(a_{1}, b_{1}), (a_{2}, b_{2}),... (a_{m_{i}+1}, b_{m_{i}+1})$. The multiple integral $I_{m_{j}+1} ^{j}$ acts with respect to the variables $(c_{1},d_{1}), (c_{2}, d_{2}),..., (c_{m_{j}+1}, d_{m_{j}+1})$. Due to the independence of the components $W^{(i)}$, $i=1,..,d$, the term $\Delta _{i,j} (x-W _{a_{k}, d_{l}})$, although stochastic,  can be viewed as a deterministic integrand for both multiple Wiener-It\^o integrals $I_{m_{i}+1}^{i}$ and $I_{m_{j}+1} ^{j}$. Also the integrals $I_{m_{i}+1}^{i}$ and $I_{m_{j}+1} ^{j}$ plays the role of a deterministic integrand one for the other.

Let us compute the Fourier transform of $\xi_{i,j}$, for $1\leq i,j\leq d$, $i\not=j$. We recall that the Fourier transform is taken with respect to the variable $x\in \mathbb{R} ^{d}$.

\begin{eqnarray}
\widehat{\xi } _{i,j}(x)&=& \sum_{m_{i}, m_{j}\geq 0} \frac{1}{m_{i}+ 1} \frac{1}{m_{j}+1} a_{m_{i}}^{\hat{x}_{i}}(a_{k}, d_{l})a_{m_{j}}^{\hat{x}_{j}}(a_{k}, d_{l})\nonumber \\
&& I_{m_{i}+1}^{i} I_{m_{j}+1} ^{j} \left( \sum_{k=0} ^{m_{i}+1} \sum_{l=0} ^{m_{j}+1} \Delta _{i,j}(\hat{x}-W_{a_{k} , d_{l}}) \right.\nonumber\\
&& \left. 1_{[0,a_{k}] \times [0, d_{l} ] } ^{\otimes m_{i} } \left( (a_{1}, b_{1}), ..., \overline{(a_{k}, b_{k})},..., (a_{m_{i}+1} , b_{m_{i}+1}) \right)1_{[0, a_{k}]}(c_{l}) 1_{[0, d_{l}]} (b_{k})\right. \nonumber\\
&& \left. 1_{[0,a_{k}] \times [0, d_{l} ] } ^{\otimes m_{i} }\left( (c_{1}, d_{1}), ..., \overline{(c_{l}, d_{l})},...,(c_{m_{j}+1}, d_{m_{j}+1}) \right) \right) \label{n1}
\end{eqnarray}
where $a_{m_{i}}^{\hat{x}_{i}}$ denotes the Fourier transform of the function $x\to a_{m_{i}}^{x_{i}}$ and $ \Delta _{i,j}(\hat{x}-W_{a_{k} , d_{l}})$ means the Fourier transform of the function
$$(x_{1}, .., \overline{x}_{i},.., \overline{x}_{j},.., x_{d} )\in \mathbb{R} ^{d-2} \to \Delta _{i,j}(x-W_{a_{k} , d_{l}}).$$
We have from Lemma \ref{ahat}
\begin{equation}\label{n2}
a_{m_{i}}^{\hat{x}_{i}}(a_{k}, d_{l}) =\frac{(-1) ^{m_{i}} x_{i} ^{m_{i}}}{m_{i}!} e^{-\frac{x_{i}^{2} a_{k}d_{l}}{2}}
\end{equation}
and clearly
\begin{equation}\label{n3}
 \Delta _{i,j}(\hat{x}-W_{a_{k} , d_{l}})=\prod _{k=1; k\not= i,j} ^{d} e^{-ix_{k} W^{(k)}_{a_{k}, d_{l}}}.
\end{equation}
Then, by combining (\ref{n1}), (\ref{n2}), (\ref{n3})
\begin{eqnarray*}
\widehat{\xi } _{i,j}(x)&=& \sum_{m_{i}, m_{j}\geq 0} \frac{1}{m_{i}+ 1} \frac{1}{m_{j}+1}\frac{(-i) ^{m_{i}+ m_{j}}}{m_{i}! m_{j}! }x_{i}^{m_{i}}x_{j}^{m_{j}}I_{m_{i}+1}^{i} I_{m_{j}+1} ^{j} \\
&&\sum_{k=1} ^{m_{i}+1} \sum_{l=1} ^{m_{j}+1} e^{-\frac{x_{i}^{2} a_{k}d_{l}}{2}}e^{-\frac{x_{j}^{2} a_{k}d_{l}}{2}}\left( \prod _{k=1; k\not= i,j} ^{d} e^{-ix_{k} W^{(k)}_{a_{k}, d_{l}}}\right) \\
&& 1_{[0,a_{k}] \times [0, d_{l} ] } ^{\otimes m_{i} } \left( (a_{1}, b_{1}), ..., \overline{(a_{k}, b_{k})},..., (a_{m_{i}+1} , b_{m_{i}+1}) \right)1_{[0, a_{k}]}(c_{l})\\
&&  1_{[0,a_{k}] \times [0, d_{l} ] } ^{\otimes m_{j} }\left( (c_{1}, d_{1}), ..., \overline{(c_{l}, d_{l})},...,(c_{m_{j}+1}, d_{m_{j}+1}) \right)1_{[0,d_{l}]} (b_{k})\\
&=&\sum_{m_{i}, m_{j}\geq 0} \frac{1}{m_{i}+ 1} \frac{1}{m_{j}+1}\frac{(-i) ^{m_{i}+ m_{j}}}{m_{i}! m_{j}! }x_{i}^{m_{i}}x_{j}^{m_{j}}I_{m_{i}+1}^{i} I_{m_{j}+1} ^{j} \\
&&\sum_{k=1} ^{m_{i}+1} \sum_{l=1} ^{m_{j}+1} e^{-\frac{x_{i}^{2} a_{k}d_{l}}{2}}e^{-\frac{x_{j}^{2} a_{k}d_{l}}{2}}\left( \prod _{k=1; k\not= i,j} ^{d} e^{-ix_{k} W^{(k)}_{a_{k}, d_{l}}}\right) \\
&&1_{[0, a_{k}] } ^{\otimes m_{i}}(a_{1},.., \overline{a_{k}},.., a_{m_{i}+1}) 1_{[0, a_{k}] } ^{\otimes m_{j}+1}(c_{1},...., c_{m_{j}+1})\\
&& 1_{[0,d_{l}] }^{m_{i} +1} (b_{1},..., b_{m_{i}+1}) 1_{[0,d_{l}] }^{m_{j} }(d_{1},.., \overline{d_{l}},.., d_{m_{j}+1})
\end{eqnarray*}
and taking the square mean, we will obtain
\begin{eqnarray*}
\mathbf{E}\left| \widehat{\xi } _{i,j}(x)\right| ^{2} &=&  \sum_{m_{i}, m_{j}\geq 0} \frac{1}{(m_{i}+ 1)!} \frac{1}{(m_{j}+1)!}x_{i}^{2m_{i}}x_{j}^{2m_{j}} \sum_{k=1} ^{m_{i}+1} \sum_{l=1} ^{m_{j}+1}e ^{-x_{i} ^{2}a_{k}d_{l} }e ^{-x_{j} ^{2}a_{k}d_{l} }\\
&& \int_{[0,1] ^{m_{i}+1} }da_{1}...da_{m_{i}+1} \int_{[0,1] ^{m_{i}+1} }db_{1}...db_{m_{i}+1} \\
&&\int_{[0,1] ^{m_{j}+1} } dc_{1}...dc_{m_{j+1} } \int_{[0,1] ^{m_{i}+1} }dd_{1}...dd_{m_{j}+1} \\
&&1_{[0, a_{k} ]} ^{\otimes m_{i}}(a_{1},.., \hat{a_{k}}, ..., a_{m_{i}+1} ) 1_{[0, d_{l}]} ^{\otimes m_{i}+1} (b_{1},...., b_{m_{i}+1})\\
&&
1_{[0, a_{k} ]} ^{\otimes m_{j}+1}(c_{1},...., c_{m_{j}+1} )1_{[0, d_{l}]} ^{\otimes m_{j}} (d_{1},.., \hat{d_{l}},.., d_{m_{i}+1})\\
&=&\sum_{m_{i}, m_{j}\geq 0} \frac{1}{(m_{i}+ 1)!} \frac{1}{(m_{j}+1)!}x_{i}^{2m_{i}}x_{j}^{2m_{j}} \\
&&\sum_{k=1} ^{m_{i}+1} \sum_{l=1} ^{m_{j}+1}\int_{0}^{1} da_{k} \int_{0}^{1}dd_{l} (a_{k}d_{l} ) ^{m_{i}+ m_{j}+1}e ^{-x_{i} ^{2}a_{k}d_{l} }e ^{-x_{j} ^{2}a_{k}d_{l} }\\
&=&\sum_{m_{i}, m_{j}\geq 0} \frac{1}{(m_{i})!} \frac{1}{(m_{j})!}x_{i}^{2m_{i}}x_{j}^{2m_{j}}  \int_{0}^{1}da \int_{0}^{1}db e^{-(x_{i}^{2}+ x_{j}^{2})ab} (ab)^{m_{i}+ m_{j}+1} .
\end{eqnarray*}
Since
\begin{equation*}
\sum_{m_{i}\geq 0} \frac{1}{m_{i}!} (x_{i}^{2} ab) ^{m_{i}} = e^{x_{i} ^{2}ab}
\end{equation*}
we get
\begin{equation*}
\mathbf{E}\left| \widehat{\xi } _{i,j}(x)\right| ^{2}=\int_{0}^{1}ada \int_{0}^{1}b db  e^{-(x_{i}^{2}+ x_{j}^{2})ab}e^{x_{i} ^{2}ab}e^{x_{j} ^{2}ab}=\frac{1}{4}.
\end{equation*}
If $i=j$ it follows from the proof of Theorem 3 that $
\mathbf{E}\left| \widehat{\xi } _{i,i}(x)\right| ^{2}=\frac{1}{4}$ for every $i=1,..,d$. As a consequence
\begin{equation*}
\int_{\mathbb{R}^{d}} \frac{1}{(1+ \vert x\vert ^{2})^{r}} \mathbf{E}\left| \widehat{\xi } (x)\right| ^{2} dx= \frac{d^{2}}{4} \int_{\mathbb{R} ^{d}} \frac{1}{(1+ \vert x\vert ^{2})^{r}}dx
\end{equation*}
is finite if and only if $r>\frac{d}{2}$. \qed

\begin{remark}
The 2D-stochastic current keeps the regularity of the one-dimensional counterpart (see \cite{FT}).
\end{remark}

 \begin{remark}
 Let $D\left(  \mathbb{R}^{d}\right)  $ be the space of smooth compact support
functions on $\mathbb{R}^{d}$ and let $D^{\prime}\left(  \mathbb{R}%
^{d}\right)  $ be its dual, the space of distributions, endowed with the usual
topologies. We denote by $\left\langle S,\varphi\right\rangle $ the dual
pairing between $S\in D^{\prime}\left(  \mathbb{R}^{d}\right)  $ and
$\varphi\in D\left(  \mathbb{R}^{d}\right)  $. Then the quantity $\langle T_{d}, \varphi\rangle $ because a well-defined random variable for every $\varphi \in D^{\prime}\left(  \mathbb{R}^{d}\right)$. Moreover, the functional $\Delta (x-W_{u_{1}, v_{2}})$ can be interpreted as a random distribution as described in Section 3 of \cite{FT}.
\end{remark}

\end{document}